\documentclass[a4paper,10pt]{article}

\usepackage{geometry}
\usepackage{epsfig}
\usepackage{amsmath}
\usepackage{amssymb}
\usepackage{amsthm}
\usepackage{mathrsfs}
\usepackage{color}
\usepackage{amscd}
\usepackage{euscript}

\theoremstyle{plain} 
\newtheorem{thm}{Theorem}[section] 
\newtheorem{cor}[thm]{Corollary} 
\newtheorem{lem}[thm]{Lemma} 
\newtheorem{prop}[thm]{Proposition} 
\theoremstyle{definition} 
\newtheorem{defn}{Definition}[section] 
\newtheorem{oss}{Remark}

\newtheorem{ex}{Example}

\DeclareFontEncoding{LS1}{}{}
\DeclareFontSubstitution{LS1}{stix}{m}{n}
\DeclareSymbolFont{symbols2}{LS1}{stixfrak} {m} {n}
\DeclareMathSymbol{\operp}{\mathbin}{symbols2}{"A8}

\usepackage{amsfonts}
\usepackage{mathtools}

\newcommand{\hslashslash}{%
  \raisebox{.9ex}{%
    \scalebox{.7}{%
      \rotatebox[origin=c]{18}{$-$}%
    }%
  }%
}
\newcommand{\fslash}{%
  {%
   \vphantom{f}%
   \ooalign{\kern.05em\smash{\hslashslash}\hidewidth\cr$f$\cr}%
   \kern.05em
  }%
}

\usepackage[backend=bibtex,style=numeric]{biblatex}
\usepackage{ stmaryrd }

\usepackage{tikz-cd}

\title{{\bf Slice regular holomorphic Cliffordian functions of order $k$ }}

\author{Giulio Binosi \\
\small Dipartimento di Matematica, Universit\`a di Trento\\ 
\small Via Sommarive 14, I-38123 Povo Trento, Italy\\
\small giulio.binosi@unitn.it
}

\date{  }

\begin{document}

\maketitle


\begin{abstract}
Holomorphic Cliffordian functions of order $k$ are functions in the kernel of the differential operator $\overline{\partial}\Delta^k$. 
When $\overline{\partial}\Delta^k$ is applied to functions defined in the paravector space of some Clifford Algebra $\mathbb{R}_m$ with an odd number of imaginary units, the Fueter-Sce construction establishes a critical index $k=\frac{m-1}{2}$ (sometimes called Sce exponent) for which the class of slice regular functions is contained in the one of holomorphic Cliffordian functions of order $\frac{m-1}{2}$.
In this paper, we analyze the case $k<\frac{m-1}{2}$ and find that the polynomials of degree at most $2k$ are the only slice regular holomorphic Cliffordian functions of order $k$.

\end{abstract}


\section{Introduction}\label{sec:introduction}
In the last century, various attempts have been made to generalize the theory of holomorphic functions on the complex plane to higher dimensional algebras. Among the pioneering successes in this pursuit was the work of R. Fueter \cite{Fueter}, who introduced a comprehensive hypercomplex theory by generalizing the Cauchy-Riemann operator to the quaternionic case. In his formulation, the quaternionic regular (or Fueter-regular) functions were defined as null solutions of the extended Cauchy-Riemann operator $\overline{\partial}_{CRF}=\frac{1}{2}(\frac{\partial}{\partial\alpha}+i\frac{\partial}{\partial\beta}+j\frac{\partial}{\partial\gamma}+k\frac{\partial}{\partial\delta})$. The theory of Fueter-regular function was then generalized to any Clifford algebra $\mathbb{R}_m$, with the classical theory of Clifford analysis \cite{CliffordAnalysis}. In this context, regular functions, known as monogenic functions, were defined as those functions in the kernel of the, so called, Dirac operator\footnote{Sometimes the factor $\frac{1}{2}$ is omitted.} $\overline{\partial}\coloneqq\frac{1}{2}(\partial_{x_0}+\sum_{j=1}^m\partial_{x_j})$.

However, a significant challenge arose in this theory: in contrast with the complex case, the product of monogenic functions was not guaranteed to be monogenic. Hence, the production of examples of monogenic functions became an issue. In order to solve this problem, Fueter devised a method to generate quaternionic Fueter regular functions in two steps, starting from holomorphic functions defined on an open subset of the complex upper half plane. Given a holomorphic function $F:D\subset\mathbb{C}^+\to\mathbb{C}$, $F(\alpha+i\beta)=F_0(\alpha,\beta)+i F_1(\alpha,\beta)$, the first step of Fueter's machinery produces the quaternionic valued function $f=\mathcal{I}(F)$, defined for any $q\in\mathbb{H}$ by
\begin{equation}
\label{struttura slice functions}
    f(q)=f(q_0+\operatorname{Im}(q))=F_0(q_0,|\operatorname{Im}(q)|)+\frac{\operatorname{Im}(q)}{|\operatorname{Im}(q)|}F_1(q_0,|\operatorname{Im}(q)|).
\end{equation}
The second step of Fueter's construction consists in applying the four-dimensional Laplacian to the induced quaternionic function yielding a monogenic function $\Delta_4 f$. This construction was then extended by Sce \cite{sce} to Clifford algebras $\mathbb{R}_m$, with $m$ odd, replacing $\Delta_4$ with $\Delta_{m+1}^{\frac{m-1}{2}}$ and finally completed by Qian \cite{qian} in the case when $m$ is even, using techniques of Fourier analysis, dealing with fractional powers of the Laplacian. 

The class of functions generated by the first step of the Fueter-Sce construction has recently been considered, as a special case within the theory of slice regular functions  \cite{Struppa1},\cite{Struppa2}\cite{libroCaterina}. This theory was initially explored in the quaternionic case, then generalized to any real alternative $\ast$-algebra \cite{SRAA}. It includes polynomials and power series. Restrictions of $\mathbb{R}_m$-valued slice regular functions to the paravector space $\mathbb{R}^{m+1}$ form the class of slice monogenic functions, already constructed in \cite{sliceMonogenic}. 

Another function theory that overcomes the algebraic challenge encountered by monogenic functions is the theory of holomorphic Cliffordian functions, founded by Laville and Ramadanoff \cite{LavilleRamadanoff}, which are null solutions of a higher order differential operator. More precisely, for any odd number $m$, holomorphic Cliffordian functions are sufficiently differentiable function $f:\mathbb{R}^{m+1}\to\mathbb{R}_m$ in the kernel of $\overline{\partial}\Delta^{\frac{m-1}{2}}$. The class of holomorphic Cliffordian functions obviously contains the monogenic functions, since the Laplacian is a real operator. Moreover it also contains (restrictions to $\mathbb{R}^{m+1}$ of) slice regular functions, as established by the  Fueter-Sce theorem (see Theorem \ref{Fueter-Sce theorem}).
Recently, a broader class of functions has been considered: holomorphic Cliffordian functions of order $k$ \cite{FineStructure}. They are functions in the kernel of $\overline{\partial}\Delta^k$, for $0\leq k\leq \frac{m-1}{2}$. This paper investigates the kernel of this operator restricted to slice regular functions. It is found (Theorem \ref{Main theorem}) that, for any $k$ less than the critical index $\frac{m-1}{2}$ (often referred to as Sce exponent), polynomials of degree at most $2k$ are the only slice regular holomorphic Cliffordian functions of order $k$. Conversely, for $k\geq\frac{m-1}{2}$, every slice regular function is holomorphic Cliffordian of order $k$ again by the Fueter-Sce theorem.

\section{Preliminaries}
\subsection{Monogenic functions}
Let $m\in\mathbb{N}$ and let $\{e_0,e_1,...,e_m\}$ be an orthonormal frame of $\mathbb{R}^{m+1}$. Let us define the following product rule
\begin{equation}
\label{eq product clifford}
    \begin{split}
        &e_i\cdot e_0=e_0\cdot e_i=e_i,\qquad\forall i=1,...,m\\
        &e_i\cdot e_j+e_j\cdot e_i=-2\delta_{ij},\qquad\forall i,j=1,...,m.
    \end{split}
\end{equation}
The Clifford algebra $\mathbb{R}_m$ is the vector space of dimension $2^m$ generated by $\{e_0,e_1,...,e_m,e_1\cdot e_2,\dots,e_1\cdot\dots\cdot e_m\}$, endowed with the product defined in \eqref{eq product clifford}, extended by associativity and linearity to the whole algebra. Any element $x\in\mathbb{R}_m$ can be uniquely written as $x=\sum_{A\subset\{1,...,m\}}x_Ae_A$, where $x_A\in\mathbb{R}$ and if $A=\{i_1,...,i_k\}$, with $1\leq i_1<\dots<i_k$, $e_A\coloneqq e_{i_1}\cdot\dots\cdot e_{i_k}$. We write $e_0=e_\emptyset=1$, the unity of the algebra and identify $\mathbb{R}$ as the subalgebra generated by $1$.
 We call $x\in\mathbb{R}_m$ a paravector if $x=x_0+\sum_{|A|=1}x_Ae_A=x_0+\sum_{j=1}^mx_je_j$. The subspace of $\mathbb{R}_m$ consisting of paravectors is isomorphic to $\mathbb{R}^{m+1}$ by the isomorphism $x_0+\sum_{j=1}^mx_je_j\mapsto (x_0,x_1,\dots,x_m)$. For paravectors we will sometimes use the notation $x=x_0+\operatorname{Im}(x)\in\mathbb{R}\times\mathbb{R}^m$. 

Let us now give the definition of monogenic function.
\begin{defn}
    Let $\Omega\subset\mathbb{R}^{m+1}$ be an open set and let $\overline{\partial}\coloneqq\frac{1}{2}(\partial_{x_0}+\sum_{j=1}^k\partial_{x_i})$ be the Dirac operator of $\mathbb{R}^{m+1}$. A function $f:\Omega\to\mathbb{R}_m$ of class $\mathcal{C}^1$ is called monogenic if $\overline{\partial}f=0$. We will denote by $\mathcal{M}(\Omega)$ the set of monogenic functions with domain $\Omega$.
\end{defn}

\subsection{Slice regular functions}
We recall some definitions from \cite{SRAA} and we adapt them to the case $A=\mathbb{R}_m$.
Consider the Clifford conjugation $x\mapsto x^c$, which is the unique involution of $\mathbb{R}_m$ such that $e_i^c=-e_i$, for any $i=1,\dots,m$. This makes $\mathbb{R}_m$ a $^\ast$-algebra. Define two maps $t,n:\mathbb{R}_m\to\mathbb{R}_m$ by
\begin{equation*}
t(x)\coloneqq x+x^c,\qquad n(x)=xx^c.
\end{equation*}
Let us define the "sphere" of imaginary units $\mathbb{S}_{\mathbb{R}_m}\coloneqq\{x\in\mathbb{R}_m:t(x)=0,n(x)=1\}$ and the quadratic cone $Q_{\mathbb{R}_m}\coloneqq\{x\in\mathbb{R}_m: t(x),n(x)\in\mathbb{R}\}$.  Equivalently, we can write
$$Q_{\mathbb{R}_m}=\bigcup_{J\in\mathbb{S}_{\mathbb{R}_m}}\mathbb{C}_J,$$ where $\mathbb{C}_J$ is the subalgebra of $\mathbb{R}_m$, generated by $e_0$ and $J$. Since $J^2=-1$, for any $J\in\mathbb{S}_{\mathbb{R}_m}$, we observe that $\mathbb{C}_J\cong\mathbb{C}$. Thus, any element $x\in Q_{\mathbb{R}_m}$ can be written as $x=\alpha+J\beta$, for some $\alpha,\beta\in\mathbb{R}$ and $J\in\mathbb{S}_{\mathbb{R}_m}$. If $x\in Q_{\mathbb{R}_m}\setminus\mathbb{R}$, the choices of $\alpha\in\mathbb{R}$, $\beta\in\mathbb{R}^+$ and $J\in\mathbb{S}_{\mathbb{R}_m}$ are unique.

Let $\mathbb{R}_m\otimes\mathbb{C}=\{a+ib:a,b\in\mathbb{R}_m\}$ be the complexification of $\mathbb{R}_m$ and equip it with the conjugation defined by $\overline{a+ib}=a-ib$, for any $a,b\in\mathbb{R}_m$.
\begin{defn}
Let $D\subset\mathbb{C}$ be an open set that is invariant with respect to conjugation, namely $\overline{z}\in D$ for every $z\in D$. A function $F:D\to\mathbb{R}_m\otimes\mathbb{C}$ is called stem function if 
\begin{equation*}
    F(\overline{z})=\overline{F(z)},\qquad\forall z\in D.
\end{equation*}
If we write $F$ in components as $F=F_0+iF_1$, $F$ is a stem function if and only if its components satisfy the following even-odd conditions w.r.t. the imaginary part of $z$:
\begin{equation*}
    F_0(\overline{z})=F_0(z),\qquad F_1(\overline{z})=-F_1(z),\qquad\forall z\in D.
\end{equation*}
\end{defn}
For given $D\subset\mathbb{C}$, we define the circularization of $D$ in $Q_{\mathbb{R}_m}$ as
\begin{equation*}
    \Omega_D\coloneqq\{\alpha+J\beta: \alpha+i\beta\in D, J\in\mathbb{S}_{\mathbb{R}_m}\}\subset Q_{\mathbb{R}_m}.
\end{equation*}
Such sets $\Omega_D$ will be called circular. We also say that $\Omega_D$ is a slice (resp. product) domain if $D$ is connected and $D\cap\mathbb{R}\neq\emptyset$ (resp. $D=D^+\cup D^-$, with $D^-=\overline{D^+}$ and $D\cap\mathbb{R}=\emptyset$). 
Now, we can give the definition of slice regular function.
\begin{defn}
    Let $F=F_0+iF_1:D\to\mathbb{R}_m\otimes\mathbb{C}$ be a stem function. Then, $F$ uniquely induces the slice function $f=\mathcal{I}(F):\Omega_D\to\mathbb{R}_m$ defined for any $x=\alpha+J\beta$ as
    \begin{equation*}
        f(x)=F_0(\alpha+i\beta)+JF_1(\alpha+i\beta).
    \end{equation*}
    We say that $f$ is slice regular if $F$ is holomorphic with respect to the complex structure defined by left multiplication by $i$. We will denote by $\mathcal{S}(\Omega_D)$ and $\mathcal{S}\mathcal{R}(\Omega_D)$ respectively, the set of slice and slice regular functions on $\Omega_D$.
    
    \begin{oss}
    Since the paravector subspace $\mathbb{R}^{m+1}$ is always contained in $Q_{\mathbb{R}_m}$, we can restrict any function $f\in\mathcal{S}(\Omega_D)$ to $\Omega=\Omega_D\cap\mathbb{R}^{m+1}$. Thanks to representation formula (\cite[Proposition 6]{SRAA}), the restriction to $\Omega$ uniquely determines a slice function on $\Omega_D$. By this, one can also consider the sets $\mathcal{S}(\Omega)=\{f|_\Omega:f\in\mathcal{S}(\Omega_D)\}$, $\mathcal{SR}(\Omega)=\{f|_\Omega:f\in\mathcal{SR}(\Omega_D)\}$.
\end{oss}
    
    If $F:D\to\mathbb{R}_m$ is a $\mathbb{R}_m$-valued stem function (namely, if $F_1=0$), we will call the induced slice function $\mathcal{I}(F)$ circular. Circular slice functions are characterized by the fact that they are constant above spheres $\mathbb{S}_{\alpha,\beta}=\{\alpha+J\beta: J\in\mathbb{S}_{\mathbb{R}_m}\}$.
\end{defn}

Given a slice function $f:\Omega_D\to\mathbb{R}_m$ we can define its spherical derivative $f'_s:\Omega_D\setminus\mathbb{R}\to\mathbb{R}_m$, as
\begin{equation*}
    f'_s(x)=\left[2\operatorname{Im}(x)\right]^{-1}(f(x)-f(\overline{x})).
\end{equation*}
Note that if $f=\mathcal{I}(F_0+iF_1)$, $f'_s$ is a circular slice function induced by the $\mathbb{R}_m$-valued stem function $F'_s\colon D\setminus\mathbb{R}\to\mathbb{R}_m$, $F'_s(\alpha+i\beta)\coloneqq F_1(\alpha+i\beta)/\beta$. Since $F_1$ is an odd function w.r.t. $\beta$, $F'_s$ is even.
If $f\in\mathcal{SR}(\Omega_D)$ or, equivalently, if $F$ is holomorphic, then $f'_s$ continuously extends to $\Omega_D\cap\mathbb{R}$ and $F'_s$ to $D\cap\mathbb{R}$.
Thus, we have that $f'_s(\alpha+J\beta)=F'_s(\alpha+i\beta)$, for any $\alpha+i\beta\in D$ and any $J\in\mathbb{S}_{\mathbb{R}_m}$.


The spherical derivative is closely related to the Dirac operator when they act on slice regular functions. The next result derives from \cite[Corollary 3.3.3]{HarmonicityPerotti}, taking into account the slightly different definition of $\overline{\partial}$ given in the present work. 
\begin{prop}
    \label{Proposizione Perotti derivata sferica e Dirac operator}
    Let $f$ be (the restriction to $\Omega$ of) a slice regular function $\Omega_D \to \mathbb{R}_m$. Then it holds
    \begin{equation*}
        \overline{\partial}f=\frac{1-m}{2}f'_s.
    \end{equation*}
\end{prop}
The Fueter-Sce theorem constitutes a bridge between the theories of slice regular and monogenic functions. Originally, it was conceived for the subclass of slice preserving slice regular functions. In modern language, it can be stated as follows (see, for example \cite{FueterandFcalculus}).
\begin{thm}[Fueter-Sce theorem]
    \label{Fueter-Sce theorem}
    Let $m$ be odd and let $f$ be (the restriction to $\Omega$ of) a slice regular function $\Omega_D \to \mathbb{R}_m$. Then $\Delta^{\frac{m-1}{2}}f$ is axially-monogenic, i.e. $\Delta^{\frac{m-1}{2}}f\in\mathcal{S}(\Omega)\cap\mathcal{M}(\Omega)=:\mathcal{A}\mathcal{M}(\Omega)$. Namely, we have the Fueter-Sce map
    \begin{equation*}
            \Delta^{\frac{m-1}{2}}:\mathcal{S}\mathcal{R}(\Omega)\to\mathcal{A}\mathcal{M}(\Omega).
    \end{equation*}
\end{thm}

\subsection{Holomorphic Cliffordian functions}
\begin{defn}
    Let $\Omega\subset\mathbb{R}^{m+1}$ be an open set. A function $f:\Omega\to\mathbb{R}_m$ of class $\mathcal{C}^{2k+1}$ is called holomorphic Cliffordian of order $k$ if $\overline{\partial}\Delta_{m+1}^kf=0$ or, equivalently, if $\Delta_{m+1}^kf$ is monogenic. Holomorphic Cliffordian functions of order $\frac{m-1}{2}$ will be simply called holomorphic Cliffordian functions.
\end{defn}
Observe that monogenic functions are holomorphic Cliffordian of any order, since $\overline{\partial}\Delta_{m+1}^kf=\Delta_{m+1}^k\overline{\partial}f=0$, for any $f\in\mathcal{M}(\Omega)$, while by the Fueter-Sce theorem slice regular functions are holomorphic Cliffordian of order $k\geq\frac{m-1}{2}$.

\section{Main result}
From here on we assume $m$ odd and we denote with $\gamma_m\coloneqq\frac{m-1}{2}$ the Sce exponent.

\begin{defn}
For every $k<\gamma_m$, let us define $\mathfrak{F}_k$ such as the operator that acts on a slice regular function $f\in\mathcal{S}\mathcal{R}(\Omega_D)$ as $\overline\partial\Delta^k_{m+1}$ applied to the restriction of $f$ to $\Omega$. Namely, 
\begin{equation*}
\mathfrak{F}_k\colon\mathcal{S}\mathcal{R}(\Omega_D)\ni f\mapsto\overline\partial\Delta^k_{m+1}(f|_\Omega)\in\ker\Delta^{\gamma_m-k}.
\end{equation*}
\end{defn}

Let us set $A\coloneqq\mathbb{R}_m$ and let $A[x]$ denote the algebra of polynomials with coefficients in $A$. It is known that $A[x]$ can be naturally included in $\mathcal{S}\mathcal{R}(Q_A)$, indeed we can associate to each element $P(x)=\sum_{n=0}^dx^na_n\in A[x]$ the slice regular function $P|_{Q_A}:Q_A\to A,x\mapsto\sum_{n=0}^dx^na_n$. The inclusion $A[x]\hookrightarrow\mathcal{S}\mathcal{R}(Q_A)$ is an injective algebra homomorphism and the elements of the image are called slice regular polynomials. For simplicity, we will identify $P(x)$ and $P|_{Q_A}$, in order to assume $A[x]\subset\mathcal{S}\mathcal{R}(Q_A)$. Let us also denote with $A_d[x]$ the set of (slice regular) polynomials of degree less than or equal to $d$.
\begin{thm}
\label{Main theorem}
Let $f\in\mathcal{S}\mathcal{R}(\Omega_D)$ and let $k<\gamma_m=\frac{m-1}{2}$. Then (the restriction to $\Omega$ of) $f$ is holomorphic Cliffordian of order $k$ if and only if $f$ is a polynomial of degree at most $2k$. In other words,
\begin{equation*}
    \ker\mathfrak{F}_k=A_{2k}[x],\qquad\forall k<\gamma_m.
\end{equation*}
\end{thm}
\begin{oss}
   The Fueter-Sce theorem asserts that $\ker\mathfrak{F}_k=\mathcal{S}\mathcal{R}(\Omega_D)$, for any $k\geq\gamma_m$. Hence, we can consider the following chain of inclusions, that ends with the Sce exponent $\gamma_m=\frac{m-1}{2}$:
    \begin{equation*}
        \ker\mathfrak{F}_0=A\subset \ker\mathfrak{F}_1=A_2[x]\subset\dots\subset\ker\mathfrak{F}_{\frac{m-3}{2}}=A_{m-3}[x]\subset\ker\mathfrak{F}_{\gamma_m}=\mathcal{S}\mathcal{R}(\Omega_D).
    \end{equation*}
\end{oss}

If we restrict to the special case $k = 0$, we find the next result (already proven in \cite[Corollary 3.3.3]{HarmonicityPerotti}).
\begin{cor}
Let $f \in \mathcal{SR}(\Omega_D)$ and assume the restriction of $f$ to $\Omega$ to be monogenic. Then $f$ is constant.
\end{cor}

\section{Proof of Theorem \ref{Main theorem}}
Let $F=F_0+iF_1:D\to\mathbb{R}_m$ be a holomorphic stem function. Since for fixed $\alpha\in\pi_1(D)=\{\alpha\in\mathbb{R}:(\alpha,\beta)\in D\}$ the function $F_1(\alpha,\cdot)$ is odd, by Whitney's Theorem \cite[Theorem 2]{Whitney}, there exists a holomorphic function $G\colon D\to \mathbb{R}_m$ that satisfies $\beta G(\alpha,\beta^2)=F_1(\alpha,\beta)$. In particular, we see that for any $(\alpha+i\beta)\in D$,
\begin{equation*}
F'_s(\alpha,\beta)=G(\alpha,\beta^2)
\end{equation*}
Let us denote with $\partial_2 G(x,y)\coloneqq\frac{\partial G}{\partial y}(x,y)$.
The next result is proven in \cite[Theorem 3.4.1]{HarmonicityPerotti}
\begin{prop}
\label{Proposition Perotti}
Let $f$ be (the restriction to $\Omega$ of) a slice regular function $\Omega_D\to\mathbb{R}_m$ and let $\Delta_{m+1}$ be the Laplacian of $\mathbb{R}^{m+1}$. Then for every $k=1,...,[\frac{m-1}{2}]$ and any $x\in\Omega$, it holds
    \begin{equation}
    \label{Equazione Perotti}
        \Delta_{m+1}^kf'_s(x)=2^k(m-3)\cdot\dots\cdot(m-2k-1)\partial_2^kG(\operatorname{Re}(x),|\operatorname{Im}(x)|^2).
    \end{equation}
\end{prop}

The goal of the following Proposition is to express \eqref{Equazione Perotti} in term of derivatives of $F'_s$.
\begin{prop}
\label{Prop Derivate G in funzione di derivata sferica}
    Let $f$ be (the restriction to $\Omega$ of) a slice function $\Omega_D\to\mathbb{R}_m$ and let $G$ as before. Then, for any $k\in\mathbb{N}$, it holds
    \begin{equation}
    \label{Relazione G e f's}
            \partial_2^kG(\alpha,\beta^2)=2^{-k}\sum_{\ell=1}^ka_\ell^{(k)}\beta^{\ell-2k}\partial_2^\ell F'_s(\alpha,\beta),
    \end{equation}
    with $a_\ell ^{(k)}\coloneqq \dfrac{(-2)^{\ell -k}(2k-\ell -1)!}{(\ell -1)!(k-\ell )!}$. Moreover, the coefficients $a_\ell ^{(k)}$ are characterized by
    \begin{equation}
        \label{recursive relation coefficient}
        \left\{\begin{array}{l}
             a_\ell ^{(k+1)}=a_{\ell -1}^{(k)}-(2k-\ell )a_\ell ^{(k)},\qquad\forall \ell =1,...,k  \\\\
             a_k^{(k)}=1, \ a_0^{(k)}=0. 
        \end{array}\right.
    \end{equation}
    Finally, if $f$ is (the restriction of) a slice regular function, $\eqref{Equazione Perotti}$ becomes
    \begin{equation}
    \label{eq laplaciano con derivata sferica}
        \Delta^k_{m+1}f'_s(x)=(m-3)\dots(m-2k-1)\sum_{\ell =1}^ka_\ell ^{(k)}|\operatorname{Im}(x)|^{\ell -2k}\partial_2^\ell F'_s(\operatorname{Re}(x),|\operatorname{Im}(x)|).
    \end{equation} 
\end{prop}

\begin{proof}
First, let us prove that $a_\ell ^{(k)}$ satisfy \eqref{recursive relation coefficient}.
\begin{equation*}
    \begin{split}
        a_{\ell -1}^{(k)}-(2k-\ell )a_k^{(k)}&=\frac{(-2)^{\ell -k-1}(2k-\ell )!}{(\ell -2)!(k-\ell +1)!}-\frac{(-2)^{\ell -k}(2k-\ell )!}{(\ell -1)!(k-\ell )!}\\
        &=\frac{(-2)^{\ell -k-1}(2k-\ell )!}{(\ell -1)!(k-\ell +1)!}\left[(\ell -1)-(-2)(k-\ell +1)\right]\\
        &=\frac{(-2)^{\ell -k-1}(2k-\ell +1)!}{(\ell -1)!(k-\ell +1)!}=a_\ell ^{(k+1)}.
    \end{split}
\end{equation*}
Now, let us prove \eqref{Relazione G e f's}. Note that, for every $ \alpha,  \gamma\in\mathbb{R}$, it holds $G( \alpha,  \gamma)=F'_s\left( \alpha,\sqrt{  \gamma}\right)$. Assume by induction that, for some $k\in\mathbb{N}$, it holds 
    \begin{equation*}
        \partial_2^kG( \alpha,  \gamma)=2^{-k}\sum_{\ell =1}^ka_\ell ^{(k)}  \gamma^{\frac{\ell -2k}{2}}\partial_2^\ell   F'_s\left( \alpha,\sqrt{  \gamma}\right),
    \end{equation*}
    so
    \begin{equation*}
        \begin{split}
            \partial_2^{k+1}G( \alpha,  \gamma)&=\frac{\partial}{\partial  \gamma}\left[2^{-k}\sum_{\ell =1}^ka_\ell ^{(k)}  \gamma^{\frac{\ell -2k}{2}}\partial_2^\ell  F'_s\left( \alpha,\sqrt{  \gamma}\right)\right]\\
            &=2^{-k}\sum_{\ell =1}^ka_\ell ^{(k)}\left[\frac{\ell -2k}{2}  \gamma^{\frac{\ell -2k-2}{2}}\partial_2^{\ell} F'_s\left( \alpha,\sqrt{  \gamma}\right)+\frac{1}{2\sqrt{  \gamma}}  \gamma^{\frac{\ell -2k}{2}}\partial_2^{\ell +1} F'_s\left( \alpha,\sqrt{  \gamma}\right)\right]\\
            &=2^{-k-1}\sum_{\ell =1}^ka_\ell ^{(k)}\left[(\ell -2k)  \gamma^{\frac{\ell -2k-2}{2}}\partial_2^{\ell} F'_s\left( \alpha,\sqrt{  \gamma}\right)+  \gamma^{\frac{\ell -2k-1}{2}}\partial_2^{\ell+1} F'_s\left( \alpha,\sqrt{  \gamma}\right)\right]
        \end{split}
    \end{equation*}
    and, by setting $\beta^2=\gamma$
    \begin{equation*}
        \begin{split}
            \partial_2^{k+1}G( \alpha,  \beta^2)&=2^{-k-1}\sum_{\ell =1}^ka_\ell ^{(k)}\left[(\ell -2k)  \beta^{\ell -2k-2}\partial_2^{\ell}F'_s( \alpha,  \beta)+  \beta^{\ell -2k-1}\partial_2^{\ell+1} F'_s( \alpha,  \beta)\right]\\
            &=2^{-k-1}\sum_{\ell =1}^ka_\ell ^{(k)}\left[(\ell -2k)  \beta^{\ell -2k-2}\partial_2^{\ell} F'_s( \alpha,  \beta)\right]+2^{-k-1}\sum_{\ell =2}^{k+1}a_{\ell -1}^{(k)}\left[  \beta^{\ell -2k-2}\partial_2^{\ell} F'_s( \alpha,  \beta)\right]\\
            &=2^{-k-1}\sum_{\ell =1}^{k+1}\left[a_\ell ^{(k)}(\ell -2k)+a_{\ell -1}^{(k)}\right]  \beta^{\ell -2k-2}\partial_2^{\ell} F'_s( \alpha,  \beta)\\
            &=2^{-(k+1)}\sum_{\ell =1}^{k+1}a_\ell ^{(k+1)}  \beta^{\ell -2(k+1)}\partial_2^\ell  F'_s( \alpha,  \beta),
        \end{split}
    \end{equation*}
    where we have used \eqref{recursive relation coefficient}. Finally \eqref{eq laplaciano con derivata sferica} follows from \eqref{Equazione Perotti} and \eqref{Relazione G e f's}, applied with $\alpha=\operatorname{Re}(z)$ and $\beta=|\operatorname{Im(z)}|$.
\end{proof}

The coefficients $a_\ell ^{(k)}$ are peculiar for producing a differential equation, whose solutions are polynomials with only even powers. First, let us recall two combinatorial tools we will need for the following proof. We define the hypergeometric function ${}_2F_1$ by
\begin{equation*}
{}_2F_1(a,b,c;z) \coloneqq\sum_{n =0}^{+\infty}\frac{(a)_n(b)_n}{(c)_n}\frac{z^n}{n!},
\end{equation*}
where $(a)_\ell \coloneqq h\cdot(h-1)\cdots(h-\ell +1)$ is the falling factorial or Pochhammer symbol. By following \cite[\S 8]{HypergeometricSummation}, we observe that if $a$ is a negative integer, then for any $b$ it holds
\begin{equation}\label{eq:hypergeometric}
{}_2F_1(a,b,2b,2)=0.
\end{equation}
See, for example, \cite{Bailey} for more details about hypergeometric series.

\begin{lem}
\label{Lemma equazione differenzial}
Suppose that $I\subset\mathbb{R}$ is an open interval and let $y\colon I\to\mathbb{R}$, with $y\in\mathcal{C}^k(I)$, satisfy the following linear homogeneous differential equation of order $k$
    \begin{equation}
    \label{equazione differenziale}
        \sum_{\ell =1}^ka_\ell ^{(k)}x^{\ell -1}y^{(\ell )}(x)=0,\qquad\forall x\in I,
    \end{equation}
    with $a_\ell ^{(k)}$ defined as in Proposition \ref{Prop Derivate G in funzione di derivata sferica}.
    Then $y$ is the polynomial $$y(x)=\sum_{\ell =0}^{k-1}c_\ell x^{2\ell },$$ for some $c_\ell \in\mathbb{R}$. 
    Suppose now that $I=I^+\cup I^-\subset\mathbb{R}$ is not necessarily connected, where $I^+$ is an open interval of $[0,+\infty)$ and $I^-=-I^+=\{-x: x\in I^+\}$. Then, any even solution on $I$ can be $\mathcal{C}^\infty$-extended to $\mathbb{R}$.

\end{lem}

\begin{proof}
Let us start with the case $I$ connected. The space of solutions of
\eqref{equazione differenziale} is a vector space of dimension $k$ and, since $\{1,x^2,...,x^{2(k-1)}\}$ is a set of linearly independent functions, it is enough to show that $x^{2h}$ is solution of \eqref{equazione differenziale} for $h=0,...,k-1$. 
Let us substitute $y=x^{2h}$ in the left-hand side of \eqref{equazione differenziale}:
\begin{equation*}
x^{2h}\sum_{\ell =1}^ka_\ell ^{(k)}(2h)_\ell =0\iff
    \sum_{\ell =1}^{k}(-2)^{\ell }\frac{(2k-\ell -1)!}{(\ell -1)!(k-\ell )!}(2h)_\ell =0.
\end{equation*}
It is easy to see that
\begin{equation*}
    \sum_{\ell =1}^{k}(-2)^{\ell }\frac{(2k-\ell -1)!}{(\ell -1)!(k-\ell )!}(2h)_\ell =-\frac{4h(2k-2)!}{(k-1)!}{}_2F_1(1-2h,1-k,2-2k;2).
\end{equation*}
and so, thanks to \eqref{eq:hypergeometric}, we conclude.

Now, suppose that $I=I^+\cup I^-$. If $0\in I^+$, then $I$ is connected and so $y(x)=\sum_{\ell =0}^{k-1}c_\ell x^{2\ell }$, for any $x\in\mathbb{R}$. Suppose that $0\notin I^+$, so $I$ is disconnected. $I^+$ and $I^-$ are open intervals of $\mathbb{R}$, so there exist $c_\ell ,d_\ell \in\mathbb{R}$ such that $y(x)=\sum_{\ell =0}^{k-1}c_\ell x^{2\ell }$ for any $x\in I^+$ and $y(x)=\sum_{\ell =0}^{k-1}d_\ell x^{2\ell }$ for any $x\in I^-$. Since $y$ is an even function, we have
\begin{equation*}
    \displaystyle\sum_{\ell =0}^{k-1}d_\ell x^{2\ell }=y(x)=y(-x)=\displaystyle\sum_{\ell =0}^{k-1}c_\ell x^{2\ell },\qquad \forall x\in I^-
\end{equation*}
and this holds if and only if $c_\ell =d_\ell $, for any $\ell =0,...,k-1$. Thus, $y(x)=\sum_{\ell =0}^{k-1}c_\ell x^{2\ell }$ for any $x\in\mathbb{R}$.
\end{proof}

\begin{oss}
    It is also possible to prove that for any $h\in\mathbb{N}$ it holds
    \begin{equation}\label{eq dimostrazione piu facile}
        \sum_{\ell =1}^{k}(-2)^{\ell }\frac{(2k-\ell -1)!}{(\ell -1)!(k-\ell )!}(2h)_\ell =(-4)^k\prod_{\ell =0}^{k-1}(h-\ell ).
    \end{equation}
    This would immediately prove that $x^{2h}$ is a solution of \eqref{equazione differenziale}, for any $h=0,\dots,k-1$. Unfortunatly, the proof of \eqref{eq dimostrazione piu facile} is long and technical.
\end{oss}

\begin{cor}
\label{Corollario intermedio}
    Let $m>2k+1$ and let $f$ be (the restriction to $\Omega$ of) a slice regular function $\mathcal{I}(F):\Omega_D \to \mathbb{R}_m$. Then 
    \begin{equation}
    \label{eq characterizing kernel power laplacian}
        \overline{\partial}\Delta_{m+1}^kf(x)=0,\quad\forall x\in\Omega\iff F'_s(\alpha,\beta)=\sum_{\ell =0}^{k-1}c_\ell (\alpha)\beta^{2\ell },\quad\forall(\alpha,\beta)\in D.
    \end{equation}
    In particular, $F'_s(\alpha,\cdot)$ can be extended to $\mathbb{R}$, for any $\alpha\in\pi_1(D)$.
\end{cor}
\begin{proof}
   By Proposition \ref{Proposizione Perotti derivata sferica e Dirac operator} and \eqref{eq laplaciano con derivata sferica} we have that
    \begin{equation*}
    0=\overline{\partial}\Delta_{m+1}^kf(x)=\Delta_{m+1}^k\overline{\partial}f(x)=-\gamma_m\Delta_{m+1}^kf'_s(x),\qquad\forall x\in\Omega
    \end{equation*}
    if and only if 
    \begin{equation*}
     0=\sum_{\ell =1}^ka_\ell ^{(k)}|\operatorname{Im}(x)|^{\ell -2k}\partial_2^\ell F'_s(\operatorname{Re}(x),|\operatorname{Im}(x)|),\qquad\forall x\in\Omega.
    \end{equation*}
    Let us distinguish two cases. If $\Omega_D$ is a slice domain, we immediately apply Lemma \ref{Lemma equazione differenzial} to $y(\beta)= F'_s(\alpha,\beta)$ for any fixed $\alpha\in\pi_1(D)$. Indeed, $ F'_s(\alpha,\cdot)$ is defined over $I_\alpha\coloneqq\{\beta\in\mathbb{R}:(\alpha,\beta)\in D\}$, which is connected and it satisfies \eqref{equazione differenziale} and so  $ F'_s(\alpha,\beta)=\sum_{\ell =0}^{k-1}c_\ell (\alpha)\beta^{2\ell }$, for any $\beta\in\mathbb{R}$. On the contrary, if $\Omega_D$ is a product domain then  $I_\alpha=I_\alpha^+\cup I_\alpha^-$, with $I_\alpha^+=I_\alpha\cap[0,+\infty)$ and $I_\alpha^-=-I_\alpha^+$ are both connected intervals. Moreover, $ F'_s(\alpha,\cdot)$ is an even function and so, thanks to Lemma \ref{Lemma equazione differenzial}, we conclude that $ F'_s(\alpha,\beta)=\sum_{\ell =0}^{k-1}c_\ell (\alpha)\beta^{2\ell }$, for any $\beta\in\mathbb{R}$.
\end{proof}

Now, we aim to partially reconstruct a slice regular function from its spherical derivative. First, we need a Lemma about entire functions.
\begin{lem}
\label{lemma funzione intera}
Let $F:\mathbb{C}\to\mathbb{C}$ be an entire function and let $u,v:\mathbb{R}^2\to\mathbb{R}$ be such that $F(\alpha+i\beta)=u(\alpha,\beta)+iv(\alpha,\beta)$, for all $(\alpha,\beta)\in\mathbb{R}^2$. Suppose that $v(\alpha,\beta)$ has the following form
    \begin{equation*}
        v(\alpha,\beta)=\sum_{\ell =0}^{k-1}c_\ell ^{(k)}(\alpha)\beta^{2\ell +1},
    \end{equation*}
    for some functions $c^{(k)}_\ell (\alpha)$. Then, the function $c_\ell ^{(k)}$ must be of the form 
    \begin{equation}
    \label{eq c ell  alpha}
        c_\ell ^{(k)}(\alpha)=(-1)^{\ell }(2k+1)\cdots(2\ell +2)\sum_{\eta=0}^{2k-2\ell -1}\frac{\alpha^\eta}{\eta!}s_{\eta+2\ell },\qquad \ell =0,...,k-1
    \end{equation}
    for some complex numbers $s_l\in\mathbb{C}$. As a consequence, $F$ is the polynomial
    \begin{equation*}
F(z)=\sum_{\ell =0}^{2k}\frac{(2k+1)!}{\ell !}s_{\ell -1}z^\ell ,
\end{equation*}
for some arbitrary real number $s_{-1}$.
\end{lem}

\begin{proof}
The components of an entire function are harmonic and, in order for $v$ to be harmonic, the functions $c_\ell ^{(k)}$ must satisfy the relations
    \begin{equation}
    \label{sistema per armonicità}
\left\{\begin{array}{l}
             c_\ell ^{''}(\alpha)=-(2\ell +2)(2\ell +3)c_{\ell +1}(\alpha),  \qquad \ell =0,...,k-2\\\\
      c_{k-1}^{''}=0.
\end{array}
\right.
    \end{equation}
    Indeed
    \begin{equation*}
        \begin{split}
\Delta_2 v(\alpha,\beta)&=\partial^2_\alpha v+\partial^2_\beta v=\sum_{\ell =0}^{k-1}c^{''}_\ell (\alpha)\beta^{2\ell +1}+\sum_{\ell =1}^{k-1}(2\ell +1)(2\ell )c_\ell (\alpha)\beta^{2\ell -1}\\
&=\sum_{\ell =0}^{k-1}c^{''}_{\ell }(\alpha)\beta^{2\ell +1}+\sum_{\ell =0}^{k-2}(2\ell +3)(2\ell +2)c_{\ell +1}(\alpha)\beta^{2\ell +1}\\
&=\sum_{\ell =0}^{k-2}\beta^{2\ell +1}(c_\ell ^{''}(\alpha)+(2\ell +2)(2\ell +3)c_{\ell +1}(\alpha))+c_{k-1}''(\alpha)\beta^{2k-1}.
\end{split}
    \end{equation*}
    So
    \begin{equation*}
       \Delta_2v(\alpha,\beta)=0\iff \left\{\begin{array}{l}
            c_\ell ^{''}=-(2\ell +2)(2\ell +3)c_{\ell +1},\quad \ell =1,...,k-2  \\\\
            c^{''}_{k-1}=0. 
       \end{array}\right. 
    \end{equation*}
     Note that the functions $c_\ell ^{(k)}$ defined in \eqref{eq c ell  alpha} satisfy \eqref{sistema per armonicità}: for every $\ell =0,...,k-2$
    \begin{equation*}
        c'_\ell (\alpha)=(-1)^\ell (2k+1)\cdots(2\ell +2)\sum_{\eta=0}^{2k-2\ell -1}\frac{\alpha^{\eta-1}}{(\eta-1)!}s_{\eta+2\ell },
\end{equation*}
\begin{equation*}
    \begin{split}
c^{''}_\ell (\alpha)&=(-1)^\ell (2k+1)\cdots(2\ell +2)\sum_{\eta=0}^{2k-2\ell -1}\frac{\alpha^{\eta-2}}{(\eta-2)!}s_{\eta+2\ell }\\
&=-(2\ell +2)(2\ell +3)\left[(-1)^{\ell +1}(2k+1)\cdots (2(\ell +1)+2)\sum_{\eta=0}^{2k-2(\ell +1)-1}\frac{\alpha^\eta}{\eta!}s_{\eta+2(\ell +1)}\right]\\
&=-(2\ell +2)(2\ell +3)c_{\ell +1}(\alpha)
\end{split}
    \end{equation*}
    and moreover $c_{k-1}(\alpha)=(-1)^{k-1}2k(2k+1)(s_{2k-2}+\alpha s_{2k-1})$ satisfies $c_{k-1}^{''}=0.$
    
    Now, let us prove that if $F(z)=\sum_{\ell =0}^{2k}\frac{(2k+1)!}{\ell !}s_{\ell -1}z^\ell $, then $v(\alpha,\beta)=\sum_{\ell=0}^{k-1}c_\ell ^{(k)}(\alpha)\beta^{2\ell +1}$. It holds
    \begin{equation*}
        \begin{split}
v(\alpha,\beta)&=\frac{1}{2i}(F(\alpha+i\beta)-F(\alpha-i\beta))=\frac{1}{2i}\sum_{t =0}^{2k}\frac{(2k+1)!}{t !}((\alpha+i\beta)^t -(\alpha-i\beta)^t )\\
&=\sum_{t =0}^{2k}\frac{(2k+1)!}{t !}s_{t -1}\sum_{\eta=0}^{\lfloor \frac{t -1}{2}\rfloor}\binom{t }{2\eta+1}\alpha^{t -2\eta-1}(-1)^\eta\beta^{2\eta+1}\\
&=\sum_{t =1}^k\dfrac{(2k+1)!}{(2t )!}s_{2t -1}\sum_{\eta=0}^{t -1}\binom{2t }{2\eta+1}\alpha^{2t -2\eta-1}(-1)^\eta\beta^{2\eta+1}+\\
&+\sum_{t =0}^{k-1}\frac{(2k+1)!}{(2t +1)!}s_{2t }\sum_{\eta=0}^t \binom{2t +1}{2\eta+1}\alpha^{2t -2\eta}(-1)^{\eta}\beta^{2\eta+1}\\
&=(2k+1)!\sum_{t =1}^k\sum_{\eta=0}^{t -1}\dfrac{(-1)^\eta s_{2t -1}\alpha^{2t -2\eta-1}\beta^{2\eta+1}}{(2\eta+1)!(2t -2\eta-1)!}+\\
&+(2k+1)!\sum_{t =0}^{k-1}\sum_{\eta=0}^{t }\dfrac{(-1)^\eta s_{2t }\alpha^{2t -2\eta}\beta^{2\eta+1}}{(2\eta+1)!(2t -2\eta)!}.
\end{split}
\end{equation*}
We can handle those two double sums for our purpose:
\begin{equation*}
    \begin{split}
&\sum_{t =1}^k\sum_{\eta=0}^{t -1}\dfrac{(-1)^\eta s_{2t -1}\alpha^{2t -2\eta-1}\beta^{2\eta+1}}{(2\eta+1)!(2t -2\eta-1)!}=\sum_{t =0}^{k-1}\sum_{\eta=0}^t  \dfrac{(-1)^\eta s_{2t +1}\alpha^{2t -2\eta+1}\beta^{2\eta+1}}{(2\eta+1)!(2t -2\eta+1)!}\\
=&\sum_{t =0}^{k-1}\sum_{\eta=0}^{k-t -1} \dfrac{(-1)^t s_{2\eta+2t +1}\alpha^{2\eta+1}\beta^{2t +1}}{(2t +1)!(2\eta+1)!}
\end{split}
\end{equation*}
and similarly
\begin{equation*}
    \sum_{t =0}^{k-1}\sum_{\eta=0}^{t }\dfrac{(-1)^\eta s_{2t }\alpha^{2t -2\eta}\beta^{2\eta+1}}{(2\eta+1)!(2t -2\eta)!}=\sum_{\mu =0}^{k-1}\sum_{\nu=0}^{k-\mu -1}\dfrac{(-1)^\mu s_{2\nu+2\mu}\alpha^{2\nu}\beta^{2\mu +1}}{(2\mu +1)!(2\nu)!},
\end{equation*}
so we get
\begin{equation*}
    \begin{split}
v(\alpha,\beta)&=(2k+1)!\left[\sum_{\ell =0}^{k-1}\sum_{\eta=0}^{k-\ell -1} \dfrac{(-1)^\ell s_{2\eta+2\ell +1}\alpha^{2\eta+1}\beta^{2\ell +1}}{(2\ell +1)!(2\eta+1)!}+\sum_{\ell =0}^{k-1}\sum_{\eta=0}^{k-\ell -1}\dfrac{(-1)^\ell s_{2\eta+2\ell }\alpha^{2\eta}\beta^{2\ell +1}}{(2\ell +1)!(2\eta)!}\right]\\
&=\sum_{\ell =0}^{k-1}(-1)^\ell \dfrac{(2k+1)!}{(2\ell +1)!}\beta^{2\ell +1}\left[\sum_{\eta=0}^{k-\ell -1} \dfrac{s_{2\eta+2\ell +1}\alpha^{2\eta+1}}{(2\eta+1)!}+\sum_{\eta=0}^{k-\ell -1}\dfrac{s_{2\eta+2\ell }\alpha^{2\eta}}{(2\eta)!}\right]\\
&=\sum_{\ell =0}^{k-1}(-1)^\ell \dfrac{(2k+1)!}{(2\ell +1)!}\beta^{2\ell +1}\sum_{\eta=0}^{2k-2\ell -1}\dfrac{s_{\eta+2\ell }\alpha^{\eta}}{\eta!}=\sum_{\ell =0}^{k-1}c_\ell ^{(k)}(\alpha)\beta^{2\ell +1},
\end{split}
\end{equation*}
as desired.
\end{proof}

We are ready to prove Theorem \ref{Main theorem}.

\begin{proof}[Proof of Theorem \ref{Main theorem}]
    Let $f=\mathcal{I}(F)$, with $F=F_0+iF_1$. From Corollary \ref{Corollario intermedio} we have
    \begin{equation*}
        \overline{\partial}\Delta_{m+1}^kf(x)=0\iff F'_s(\operatorname{Re}(x),|\operatorname{Im}(x)|)=\sum_{\ell =0}^{k-1}c_\ell (\operatorname{Re}(x))|\operatorname{Im}(x)|^{2\ell },
    \end{equation*}
    thus, 
    \begin{equation*}
        F_1(\operatorname{Re}(x),|\operatorname{Im}(x)|)=|\operatorname{Im}(x)| F'_s(\operatorname{Re}(x),|\operatorname{Im}(x)|)=\sum_{\ell =0}^{k-1}c_\ell (\operatorname{Re}(x))|\operatorname{Im}(x)|^{2\ell +1}.
    \end{equation*}
    Finally, by Lemma \ref{lemma funzione intera}, $F(z)=\sum_{\ell =0}^{2k}z^\ell c_\ell $ and so $f(x)=\mathcal{I}(F)(x)=\sum_{\ell =0}^{2k}x^\ell c_\ell $.
\end{proof}

\begin{ex}
    Let us consider the slice regular function $f:Q_{\mathbb{R}_5}\to\mathbb{R}_5$, $f(x)=x^5$. Note that $\gamma_5=2$, so we expect that $f\in\ker\mathfrak{F}_2=\mathcal{S}\mathcal{R}(Q_{\mathbb{R}_5})$ and $f\notin\ker\mathfrak{F}_k$, for $k=0,1$, since $\deg(f)=5>2k$, for any $k<\gamma_5$. Indeed, let us compute $\mathfrak{F}_kf$, which is $-4\Delta^k_{5+1}$ applied to the restriction of $f'_s$ to $\mathbb{R}^{5+1}$, for $k=0,1,2$. Let $x=x_0+\sum_{i =1}^5 x_ie_i=\alpha+J\beta$,  with $\alpha=\operatorname{Re}(x)=x_0$, $\beta=|\operatorname{Im}(x)|=\sqrt{\sum_{i =1}^5x_i^2}$ and $J=\sum_{i =1}^5x_ie_i/\beta$, then
    \begin{equation*}
f(x)=\alpha^5-10\alpha^3\beta^2+5\alpha\beta^4+J\beta(5\alpha^4-10\alpha^2\beta^2+\beta^4).
    \end{equation*}
    Thus, for $\alpha,\beta\in\mathbb{R}$ and $I\in\mathbb{S}_{\mathbb{R}_5}\cap\mathbb{R}^{5+1}$, we have 
    \begin{equation*}
         \mathfrak{F}_0(f)(\alpha+I\beta)=-4f'_s(\alpha+I\beta)=-20\alpha^4+40\alpha^2\beta^2-4\beta^4\neq0;
    \end{equation*}
    \begin{equation*}
        \mathfrak{F}_1(f)(\alpha+I\beta)=-4\Delta_{5+1} f'_s(\alpha+I\beta)=160\alpha^2-32\beta^2\neq0;
    \end{equation*}
    \begin{equation*}
        \mathfrak{F}_2(f)(\alpha+I\beta)=-4\Delta_{5+1}^2 f'_s(\alpha+I\beta)=320-320=0.
    \end{equation*}
    \end{ex}

\begin{ex}
    Let us consider again the slice regular function $f(x)=x^5$, but with $f:Q_{\mathbb{R}_9}\to\mathbb{R}_9$. Now, $\gamma_9=4$. Since $\deg(f)=5$, we expect that $f\in\ker\mathfrak{F}_3=(\mathbb{R}_9)_6[x]\subsetneq \mathcal{S}\mathcal{R}(Q_{\mathbb{R}_9})=\ker\mathfrak{F}_4$, but $f\notin\ker\mathfrak{F}_k$, for $k=0,1,2$. Indeed, let us compute as before $\mathfrak{F}_kf$, which is $-8\Delta^k_{5+1}$ applied to the restriction of $f'_s$ to $\mathbb{R}^{9+1}$, for $k=0,1,2,3$. For $\alpha,\beta\in\mathbb{R}$ and $I\in\mathbb{S}_{\mathbb{R}_9}\cap\mathbb{R}^{9+1}$, we have 
    \begin{equation*}
        \mathfrak{F}_0(f)(\alpha+I\beta)=-8f'_s(\alpha+I\beta)=-40\alpha^4+80\alpha^2\beta^2-8\beta^4\neq0;
    \end{equation*}
    \begin{equation*}
        \mathfrak{F}_1(f)(\alpha+I\beta)=-8\Delta_{9+1} f'_s(\alpha+I\beta)=960\alpha^2-192\beta^2\neq0;
    \end{equation*}
    \begin{equation*}
        \mathfrak{F}_2(f)(\alpha+I\beta)=-8\Delta_{9+1}^2 f'_s(\alpha+I\beta)=1920-3456\neq0;
    \end{equation*}
    \begin{equation*}
        \mathfrak{F}_3(f)(\alpha+I\beta)=-8\Delta_{9+1}^3 f'_s(\alpha+I\beta)=0.
    \end{equation*}
\end{ex}

We appreciate the referees'  valued comments,  which improved the clarity of the paper.

\printbibliography

\end{document}